%
\documentstyle{amsppt}
\magnification=\magstep1
\NoRunningHeads
\NoBlackBoxes
\parindent=1em
\vsize=7.4in


%
\topmatter

\title
classifying tight Weyl-Heisenberg frames
\endtitle
\author
Peter G. Casazza and Ole Christensen \endauthor
\address
Department of Mathematics,
The University of Missouri,
Columbia, Missouri 65211,
USA
and
Mathematical Institute,
Building 303,
Technical University of Denmark,
2800 Lyngby,
DENMARK
\endaddress
\email
pete\@casazza.math.missouri.edu; olechr\@mat.dtu.dk
\endemail
\thanks
The first author was supported by NSF DMS 970618.  The second
author was supported by the Danish Research Council and would
also like to thank the University of Missouri for their
hospitality.
\endthanks
%
\abstract
A Weyl-Heisenberg frame for $L^{2}(R)$ is a frame consisting of
translates and modulates of a fixed function in $L^{2}(R)$,
i.e. $(E_{mb}T_{na}g)_{m,n\in Z}$, with $a,b>0$, and $g\in L^{2}(R)$.
In this paper we will
give necessary and sufficient conditions for this family to form a
tight WH-frame.  This allows us to write down explicitly all functions
$g$ so that $(E_{mb}T_{na}g)$ is an orthonormal basis for $L^{2}(R)$.  
These results give a simple direct classification of the
 alternate dual frames to Weyl-Heisenberg frames (a result due to Janssen).
\endabstract
\endtopmatter
\document
\baselineskip=15pt
\heading{1. Introduction}
\endheading
\vskip10pt

In 1946 D. Gabor \cite{7} introduced a technique for signal processing
which soon became a paradigm for the spectral analysis associated
with time-frequency methods.  His introduction
of the short-time (windowed) Fourier transform led
eventually to Wavelet theory.  In 1952,  
Duffin and Schaeffer
 \cite{5} introduced
frame theory which put Gabor's technique into the framework of a larger
model - that of a frame for a Hilbert space.

A sequence $(f_{i})_{i\in I}$ of elements of a Hilbert space $H$ is
called a {\bf frame} for $H$ if there are constants $A,B>0$ so that
$$
A\|f\|^{2} \le \sum_{i\in I}|<f,f_{i}>|^{2} \le B\|f\|^{2},\ \ \text{
for all}\ \ f\in H.  \tag 1.1
$$
The numbers $A,B$ are called the {\bf lower} (resp. {\bf upper}) frame
bounds. The frame is a {\bf tight frame} if $A=B$ and a {\bf normalized  
tight frame} if $A = B = 1$.  If $(f_{i})$ is a frame, the frame operator is defined by
$S:H\rightarrow H,\ Sf = \sum <f,f_{i}>f_{i}$.  Then $S$ is bounded
and invertible, and
$$
f = \sum_{i\in I}<S^{-1}f,f_{i}>f_{i} = \sum_{i\in I}<f,S^{-1/2}f_{i}>
S^{-1/2}f_{i}.  \tag 1.2
$$

One can check that $(f_{i})_{i\in I}$ is a normalized tight
frame if and only if $S = I$.
One interpretation of (1.2) is that for an arbitrary frame
$(f_{i})$, a
normalized tight frame equivalent to $(f_{i})_{i\in I}$ is given by 
$(S^{-1/2}f_{i})_{i\in I}$.  Recall that Two frames $(f_{i})$ and $(g_{i})$ for a Hilbert
space $H$ are {\bf equivalent} if the operator $T:H\rightarrow H$ 
given by $T(f_{i}) = g_{i}$ is a well defined function which is an
isomorphism of $H$ onto $H$.  

If we replace $f$ in (1.1) by $f_{j}$, we see that
$$
\|f_{j}\|^{4} + \sum_{i\in I-\{j\}} | <f_{j},f_{i}>|^{2} \le B\|f_{j}\|^{2}.
$$
This yields immediately the following remark:

\proclaim{Remark 1.1}
For all $j\in I$, $\|f_{j}\|^{2}\le B$. If $\|f_{j}\|^{2} = B$  
then $f_{j}\perp \text{span}_{i\not= j}f_{i}$.  In particular, if  
$(f_{i})_{i\in I}$ is a
normalized tight frame, then $\|f_{j}\| \le 1$, and
$\|f_{j}\| = 1$ if and only if
$f_{j}\perp \text{span}_{i\not= j}f_{i}$.
\endproclaim

The particular frames of interest to us will be the Weyl-Heisenberg
frames.  To define these frames, let $a,b\in R$ and define the
operators of {\bf modulation} $E_{b}$, and {\bf translation}
$T_{a}$ for functions $f\in L^{2}(R)$ by:
$$
E_{b}f(x) = e^{2{\pi}ibx}f(x),
$$
and
$$
T_{a}f(x) = f(x-a).
$$
Given $g\in L^{2}(R)$, and $a,b >0$, we say that $(g,a,b)$ {\bf
generates a WH-frame for} $L^{2}(R)$ if $(E_{mb}T_{na}g)_{m,n\in Z}$
is a frame for $L^{2}(R)$.  The function $g$ is referred to as the
{\bf mother wavelet}.  The numbers
$a,b$ are the {\bf frame parameters}, with $a$ being the {\bf shift
parameter} and $b$ the {\bf modulation parameter}.

Let {\bf PF} denote the set of functions $g\in L^{2}(R)$ for which
$(E_{mb}T_{na}g)$ has a finite upper frame bound.  It is easily seen
that if $g\in$PF, then $g$ is bounded.  Theorem 2.1 of Casazza and  
Christensen \cite{1} yields a sufficient condition for $g\in$PF:

\proclaim{Proposition (Casazza/Christensen [1])}
If
$$
\sum_{k\in Z}|\sum_{n\in Z}g(x-na)\overline{g(x-na-k/b)}| \le B \ \  
\text{a.e.},  \tag {1.3}
$$
then $g\in$PF.
\endproclaim

It can be shown that a variation of this produces a necessary
condition for $g\in$PF.  That is, if $g\in$PF then,
$$
\sum_{k\in Z}|\sum_{n\in Z}g(x-na)\overline{g(x-na-k/b)}|^{2} < \infty.
$$
This condition is clearly
not strong enough to classify PF.  Also, condition (1.3) 
for $g\in$PF is not necessary \cite{2}. 
This condition and its relationship to the
Walnut representation of the frame operator are examined in detail
by Casazza, Christensen, and Janssen \cite{2}.

In the above proposition and in the rest of the paper, we will be working
with infinite sums of the form:
$$
\sum_{n\in Z}g(x-na)\overline{g(x-na-k/b)}.
$$
So let us discuss the question of convergence here, and ignore it for
the rest of the paper.

\proclaim{Proposition 1.2}
If $g\in L^{2}(R)$ then the series
$$
\sum_{n\in Z}g(x-na)\overline{g(x-na-k/b)}
$$
converges absolutely a.e.
\endproclaim

\demo{Proof}
Since $g,T_{k/b}g\in L^{2}(R)$ we have that $g\overline{T_{k/b}g} \in L^{1}(R)$.  Also,
$$
\|g\overline{T_{k/b}g}\|_{L^{1}} = \int_{R}|g(x)T_{k/b}(x)|\ dx =
\int_{0}^{a}\sum_{n\in Z}|g(x-na)g(x-na-k/b)|\ dx < \infty.
$$
It follows that 
$$
\sum_{n\in Z}| g(x-na)\overline{g(x-na-k/b)} | < \infty\ \ \text{a.e.}
$$
\enddemo

We need a result that appeared the first time in \cite{9}:

\proclaim{WH-Frame Identity}
If $g\in L^{2}(R)$ and
$f\in L^{2}(R)$ is bounded and compactly supported, then
$$
\sum_{n\in Z}\sum_{m\in Z}|<f,E_{mb}T_{na}g>|^{2} = F_{1}(f) +  
F_{2}(f) 
$$
where
$$
F_{1}(f) = b^{-1}\int_{R}|f(x)|^{2}\sum_{n}|g(x-na)|^{2}\ dx, 
$$
$$
F_{2}(f) = b^{-1}\sum_{k\not= 0} \int_{R}\overline{f(x)}f(x-k/b)
\sum_{n}g(x-na)\overline{g(x-na-k/b)}\ dx =  
$$
$$
b^{-1}\sum_{k\ge 1} 2\text{Re}\int_{R}\overline{f(x)}f(x-k/b)
\sum_{n}g(x-na)\overline{g(x-na-k/b)}\ dx.
$$
\endproclaim

To simplify the notation a little we introduce the following auxilliary
functions:
$$
G(x) = \sum_{n\in Z}|g(x-na)|^{2},  
$$
and for all $k\in Z$,
$$
G_{k}(x) = \sum_{n\in Z}g(x-na)\overline{g(x-na-k/b)},  
$$
It follows that $G_{0} = G$, and $G_{k}$ are periodic functions on
$R$ of period a.

\heading{2.  Classifying Tight WH-Frames}
\endheading
\vskip10pt

Now we will classify the tight WH-frames both abstractly in terms
of the behavior of related families of vectors and concretely in
terms of the behavior of the functions $G_{k}(x)$.  In the next
section we will write down explicitly the functions which satisfy
the conditions of our theorem for some cases.  Some of  

We start with a  
basic fact
which will simplify parts of the theorem and will be used in section
4 to classify the alternate dual frames for a WH-frame.

\proclaim{Proposition 2.1}
Let $g,h\in L^{2}(R)$ and $a,b\in R$.
\vskip5pt
(1) $h\perp E_{mb}g$ for all $m\not= 0$ if and only if there is a
constant $C$ so that
$$
\sum_{n\in Z}h(x-n/b)\overline{g(x-n/b)} = C, \ \ \text{a.e.}.
$$
\vskip5pt
(2)  If $n\not= 0$, then $h\perp E_{mb}T_{na}g$, for all $m\in Z$ if
and only if,
$$
\sum_{k}h(x-k/b)\overline{g(x-k/b-na)} = 0, \ \ \text{a.e.}.
$$
\endproclaim

\demo{Proof}
$(2)$: We just calculate: 
$$
 <h,E_{mb}T_{na}g> = \int_{R}h(x)\overline{E_{mb}g(x-na)}\ dx =
\int_{R}h(x)\overline{g(x-na)}E_{-mb}\ dx = 
$$
$$
\int_{0}^{1/b}\sum_{k\in Z}h(x-k/b)\overline{g(x-k/b-na)}E_{-mb}\ dx.
$$
(2) follows immediately from here and (1) is just the case of n=0 in
(2).
\enddemo

We are now ready to prove the classification theorem for tight WH-frames.
The equivalence of (1) and (2) in this theorem was first done by
Janssen \cite{6}, Section 1.3.2.  Also, the equivalence of (1) and (3)
could be derived from results of Janssen \cite{11}, Theorem 3.1.  But
our proof is is simple and follows directly from the definitions.

\proclaim{Theorem 2.2}
Let $g\in L^{2}(R)$ and $a,b\in R$.  The following are equivalent:
\vskip5pt
(1)  $(E_{mb}T_{na}g)_{n,m\in Z}$ is a normalized tight Weyl-Heisenberg frame
for $L^{2}(R)$.
\vskip5pt
(2)  We have:

\ \ \ \ \ (a)  $G(x) = \sum_{n\in Z}|g(x-na)|^{2} = b$ a.e.,
\vskip5pt
\ \ \ \ \ (b)  $G_{k}(x) = \sum_{n\in  
Z}g(x-na)\overline{g(x-na-k/b)}= 0$ a.e. for all $k\not= 0$.
\vskip5pt
(3)  We have, $g\perp E_{n/a}T_{m/b}g$, for all $(n,m)\not= (0,0)$,
and $\|g\|^{2} = ab$.
\vskip5pt
(4)  $(E_{n/a}T_{m/b}g)_{n,m\in Z}$ is an orthogonal sequence in
$L^{2}(R)$ and $\|g\|^{2} = ab$.
\vskip5pt
(5)  $(E_{mb}T_{na}g)_{n,m\in Z}$ is a Weyl-Heisenberg frame for 
$L^{2}(R)$ with frame operator $S$ and $Sg = g$.
\vskip5pt
Moreover, in case the conditions are satisfied, 
$(E_{mb}T_{na}g)_{m,n\in Z}$ is an orthonormal basis for
$L^{2}(R)$ if and only if $\|g\|=1$.
\endproclaim

\demo{Proof}
$(1)\Rightarrow (2)$:  Assume $(E_{mb}T_{na}g)_{n,m\in Z}$ is a  
normalized tight
 frame for $L^{2}(R)$. For any function $f\in L^{2}(R)$ which is
bounded and 
supported on an interval of length $< 1/b$ we have for all $x\in R$ 
and all $0\not= k\in Z$ that  $\overline{f(x)}f(x-k/b) =0$.  That is,
$F_{2}(f) = 0$.  Now, by the WH-frame Identity:
$$
\|f\|^{2} = \int_{R}|f(x)|^{2}\ dx =
\sum_{n}\sum_{m}|<f,E_{mb}T_{na}g>|^{2} = F_{1}(f) + F_{2}(f)
$$
$$
= b^{-1}\int_{R}|f(x)|^{2}\sum_{n}|g(x-na)|^{2}\ dx =
b^{-1}\int_{R}|f(x)|^{2}G(x)\ dx.
$$
Since this equality holds for all bounded $f\in L^{2}(I)$, for any interval $I$
of length $<1/b$,
it follows easily that $G(x) = b$, a.e.  Hence, for all $f\in L^{2}(R)$,
$F_{1}(f) = \|f\|^{2}$.  But now, again by the WH-frame Identity,  
we have for all bounded, compactly supported $f\in L^{2}(R)$,
$$
\|f\|^{2} = F_{1}(f)+F_{2}(f) = \|f\|^{2}+F_{2}(f).
$$
That is, $F_{2}(f) = 0$, for all bounded compactly supported $f\in  
L^{2}(R)$.
  Now fix $k_{0}\ge 1$ and
let $I$ be any interval in $R$ of length $\le 1/b$.  Define a function
$f\in L^{2}(R)$ by:
$$
f(x) = e^{i\ \text{arg}G_{k_{0}}(x)},\ \ \text{for all}\ \ x\in I,
$$
and $f(x-k_{0}/b) = 1$ for all $x\in I$ and $f(x) = 0$, otherwise.  Then,
by the WH-frame Identity,
$$
0 = F_{2}(f) = b^{-1}\sum_{k\ge 1}  
2\text{Re}\int_{R}\overline{f(x)}f(x-k/b)
\sum_{n}g(x-na)\overline{g(x-na-k/b)}\ dx =
$$
$$
b^{-1}2\text{Re}\int_{R}\overline{f(x)}f(x-k_{0}/b)
G_{k_{0}}(x)\ dx = b^{-1}2\int_{I}|G_{k_{0}}(x)|\ dx
$$
It follows that $G_{k_{0}}(x) = 0$, a.e. on $I$.  Since $k_{0}$
and $I$ were arbitrary, we have $(2b)$.

$(2)\Rightarrow (1)$:  By assumption $(2b)$ and the WH-frame Identity,
we have $F_{2}(f) = 0$ for all bounded, compactly supported $f\in  
L^{2}(R)$.
  Hence, applying
assumption $(2a)$ and the WH-frame Identity:
$$
\sum_{n}\sum_{m}|<f,E_{mb}T_{na}g>|^{2} = F_{1}(f) =
b^{-1}\int_{R}|f(x)|^{2}\sum_{n}|g(x-na)|^{2}\ dx =
$$
$$
\int_{R}|f(x)|^{2}\ dx
= \|f\|^{2}.
$$
Since this equality holds on a dense subset of $L^{2}(R)$, it holds
for all $f\in L^{2}(R)$.  So $(E_{mb}T_{na}g)_{n,m\in Z}$ is a
normalized tight frame for $L^{2}(R)$.

$(2)\Leftrightarrow (3)$:  By (2) of Proposition 2.1, (2b) is  
equivalent to
$g\perp E_{n/a}T_{m/b}g$, for all $m\not= 0$.  Also, we have,
$$
\|g\|^{2} = \int_{R}|g(x)|^{2}\ dx = \int_{0}^{a}\sum_{n\in  
Z}|g(x-na)|^{2}\ dx
= \int_{0}^{a}G(x)\ dx .
$$
So the rest of the equivalence follows from (1) of Proposition 2.1.

$(3)\Leftrightarrow (4)$:  This is immediate from the observation that
for all $m,n,\ell ,k\in Z$ we have:
$$
<E_{n/a}T_{m/b}g,E_{k/a}T_{{\ell}/b}g> =  
e^{2{\pi}i\frac{n-k}{a}\frac{m}{b}}
<g,E_{\frac{k-n}{a}}T_{\frac{{\ell}-m}{b}}g>.
$$

$(1)\Leftrightarrow (5)$:  Since $S$ commutes with $E_{mb}T_{na}$,
we have that $S=I$ if and only if $Sg = g$.

For the moreover part of the theorem, we just observe that for all
$m,n\in Z$, $\|g\| = \|E_{mb}T_{na}g\|$.  Hence, if $\|g\| = 1$, then
our tight frame consists of norm 1 elements which now form an orthonormal
basis by Remark 1.1.
\enddemo

With the characterization in Theorem 2.2, we can now recover easily and
directly from the definition the basic properties of
WH-frames which formerly required some more work.
Janssen \cite{10,11} first derived (1) in Corollary 2.3 below in an  
easier fashion.
For any WH-frame $(E_{mb}T_{na}g)_{m,n\in Z}$ with frame
operator $S$, a direct computation shows that
$$
S(E_{mb}T_{na}g) = E_{mb}T_{na}Sg,\ \ \text{for all}\ \ m,n\in Z.
$$
It follows that $S^{-1/2}$ also commutes with $E_{mb}T_{na}$ and
so $(E_{mb}T_{na}S^{-1/2}g)_{m,n\in Z}$ is a normalized tight WH-frame 
which is equivalent to $(E_{mb}T_{na}g)_{m,n\in Z}$ and hence must
satisfy the conditions of Theorem 2.2.  In particular, for all
$(m,n)\not= (0,0)$ we have:
$$
<S^{-1}g,E_{n/a}T_{m/b}g> = <S^{-1/2}g,E_{n/a}T_{m/b}S^{-1/2}g> = 0.
$$

\proclaim{Corollary 2.3}
Let $(E_{mb}T_{na}g)_{m,n\in Z}$ be a WH-frame for $L^{2}(R)$.
\vskip5pt
(1)  $S^{-1}g\perp
E_{n/a}T_{m/b}g$, for all $(m,n)\not= (0,0)$ and
$$
<S^{-1}g,g> = <S^{-1/2}g,S^{-1/2}g> = \|S^{-1/2}g\|^{2} = ab \le 1.
$$
 \vskip5pt
(2)  If $ab<1$ then the WH-frame is not a Riesz basis.
\vskip5pt
(3)  If $ab=1$ then the WH-frame is a Riesz basis for $L^{2}(R)$.
\endproclaim

\demo{Proof}
$(1)$:  All of this was observed before we stated the theorem except
the last inequality $ab\le 1$ which follows immediately from Remark 1.1.

$(2)$:  If our WH-frame is exact then $(E_{mb}T_{na}S^{-1/2}g)$ is an
orthonormal basis.  Hence,
$$
1 = \|S^{-1/2}g\|^{2} = ab.
$$

$(3)$:  If $ab=1$ then $(E_{mb}T_{na}S^{-1/2}g)$ is a tight WH-frame and
hence $(E_{n/a}T_{m/b}S^{-1/2}g)$ is an orthogonal sequence by Theorem
2.2 (4).  But, $ab=1$ implies that $n/a=nb$ and $m/b=na$, so $(E_{nb}
T_{ma}S^{-1/2}g)$ is an orthogonal basis for $L^{2}(R)$.  Since  
$S^{-1/2}$
is an isomorphism, it follows that $(E_{mb}T_{na}g)$ is a Riesz basis for
$L^{2}(R)$.
\enddemo

Finally, we point out what appears to be a surprising consequence of
Theorem 2.2 which seems to indicate the existence of a relationship
between the values of a function and its Fourier Transform.  

\proclaim{Corollary 2.4}
If $g\in L^{2}(R)$ and $ab\le 1$, the following are equivalent:

(1)  The function $g$ satisfies:
$$
\sum_{n}g(x-na)\overline{g(x-na-k/b)} = 0\ \ \text{a.e. for all}\ \ 
k\not= 0.
$$
and
$$
\sum_{n}|g(x-na)|^{2} = b\ \ \text{a.e.}
$$

(2)  The Fourier transform of the function $g$ satisfies:
$$
\sum_{n}\hat{g}(x-nb)\overline{\hat{g}(x-nb-k/a)} = 0\ \ \text{a.e.  
for all}\ \
k\not= 0.
$$
and
$$
\sum_{n}|\hat{g}(x-nb)|^{2} = a\ \ \text{a.e.}
$$
\endproclaim

\demo{Proof}
Since the Fourier transform of $E_{mb}T_{na}g$ is $E_{ma}T_{nb}\hat{g}$,
 Corollary 2.4 comes from
Theorem 2.2 applied to the families $(E_{mb}T_{na}g)$ and  
$(E_{ma}T_{nb}\hat{g})$ which are normalized tight WH-frames  
together.
\enddemo

\heading{3.  The Functions $g$ Giving Tight WH-Frames}
\endheading
\vskip10pt

In this section we will give an explicit representation for the functions
giving tight WH-frames given in Theorem 2.2 (2) for the case $a=b=1$.  If $f(x,y)$ is any
function of two variables, we denote by $f_{y}(x)$ the function:
$$
f_{y}(x) = f(x,y).
$$

We start with a simple
proposition which contains the basic notions which will be used in our
characterization.

\proclaim{Proposition 3.1}
The following are equivalent:

(1)  The sequence $z = (c_{n})_{n\in Z}\in {\ell}_{2}(Z)$ is orthogonal
to all of its proper shifts and $\|z\|^{2} = b.$

(2)  The unique function $h:[0,1]\rightarrow C$ with $\hat{h}(n)=c_{n}$, 
for all $n\in Z$ has $|h(x)|^{2} = b$ a.e.

(3)  There is a measurable function $f:[0,1]\rightarrow R$ so that the
function $h$ in (2) is of the form
$$
h(x) = \sqrt{b}\ e^{2{\pi}if(x)}.
$$
\endproclaim

\demo{Proof}
$(1)\Rightarrow (2)$  If $(\hat{h}(n))=(c_{n})$ is orthogonal to all its
shifts then $h\perp e^{2{\pi}im}h$, for all $m\not= 0$.  That is, for all
$m\not= 0$ we have
$$
0 = <h,e^{2{\pi}im}h> = <|h|^{2},e^{2{\pi}im}>
$$
Hence, $|h|^{2} = C$ a.e.  Also, $\|h\|^{2} = \|(c_{n})\|^{2} = b.$   

$(2)\Rightarrow (3)$:  This is obvious.

$(3)\Rightarrow (1)$:  Given $h$ as in (3), for any $ m\in Z$ we have:
$$
<h,E_{m}h> = \int_{0}^{1}h(x)\overline{h(x)E_{m}}\ dx =
\int_{0}^{1}|h(x)|^{2}e^{-2{\pi}imx}dx =
$$
$$
\int_{0}^{1}b\ e^{-2{\pi}imx}dx = b{\delta}_{m}.
$$
So $h$ is orthogonal to $e^{2{\pi}im}h$ for all $m\not= 0$.  Hence,
$(\hat{h}(n))$ is orthogonal to all its proper shifts and
$\|(\hat{h})(n))\|^{2} = b$.
\enddemo

Now we are ready to give an explicit representation for the functions
$g$ for which $(E_{mb}T_{na}g)$ is an orthonormal basis for $L^{2}(R)$.

\proclaim{Theorem 3.2}
Let $ab=1$ and $g\in L^{2}(R)$.  The following are equivalent:

(1)  $(E_{mb}T_{na}g)$ is a orthonormal basis for $L^{2}(R)$.  

(2)  $(E_{mb}T_{na}g)$ is a normalized tight WH-frame for $L^{2}(R)$.

(3)  There is a measurable function $f:[0,1]\times [0,a)\rightarrow R$
and
$$
h(x,y) = \sqrt{b}\ e^{2{\pi}if(x,y)}
$$
so that
$$
g(y+na) = \hat{h_{y}}(n),\ \ \text{for all}\ \ y\in [0,a).
$$
\endproclaim

\demo{Proof}
This is essentially immediate from our assumptions and Theorems 2.2 and
3.1.  Since $ab=1$, we have $a=1/b$ so Theorem 2.2 (2) becomes:
\vskip5pt
\ \ \ \ \ (a)  $G(x) = \sum_{n\in Z}|g(x-na)|^{2} = b,$ a.e.,
\vskip5pt
\ \ \ \ \ (b)  $G_{k}(x) = \sum_{n\in  
Z}g(x-na)\overline{g(x-(n-k)a)}= 0$, a.e. for all $k\not= 0$.
\vskip5pt
But condition (b) is equivalent to: $z_{y} = (g(y-na))_{n\in Z}$ is  
orthogonal to
all of its proper shifts and (a) is equivalent to $\|z_{y}\|^{2} =  
b$.  By
Theorem 3.1 these conditions are equivalent to:  For each $y\in  
[0,a)$ there
is a function $f_{y}:[0,1]\rightarrow R$ and  
$h_{y}:[0,a]\rightarrow C$ with
$$
h_{y}(x) = e^{2{\pi}if_{y}(x)},
$$
and
$$
\hat{h_{y}}(n) = g(y-na).
$$
So defining $f(x,y):[0,1]\times [0,a)\rightarrow R$ and $h(x,y)$ by:
$$
f(x,y) = f_{y}(x),\ \ \ \ h(x,y) = h_{y}(x),
$$
yields the theorem modulo
the measurability conditions which are obvious.
\enddemo

\heading{4.  Alternate Dual Frames}
\endheading
\vskip10pt

If $(f_{i})_{i\in I}$ is a frame for a Hilbert space $H$, a frame  
$(h_{i})_{i\in I}$ for $H$ is called an {\bf alternate dual frame}  
or
a {\bf pseudo-dual} for
$(f_{i})_{i\in I}$ if
$$
f = \sum_{i\in I}<f,h_{i}>f_{i},\ \ \text{for all}\ \ f\in H. \tag 4.1
$$
We already know one sequence $(h_{i})_{i\in I}$ satisfying (4.1).
Namely, the sequence $(S^{-1}f_{i})_{i\in I}$.  We call $(S^{-1}f_{i})_{i\in I}$ the {\bf canonical dual} of $(f_{i})$.
If $(f_{i})_{i\in I}$ is a normalized tight frame, then $S=I$, so
the frame equals its canonical dual frame.  The converse of this clearly
holds also.  But in general, there are many alternate dual  
frames for a given frame.  For the basic properties of alternate dual frames we refer to
\cite{7,9,12}.  Now
we will use the techniques developed in section 3 to
characterize the Weyl-Heisenberg alternate dual frames for a
given Weyl-Heisenberg frame.  We need a beautiful result of  
Wexler-Raz \cite{13} (see also Janssen \cite{10,11}):

\proclaim{Theorem(Wexler-Raz [13])}
Let $g,h\in$PF.  Then $(E_{mb}T_{na}h)$ and $(E_{n/a}T_{m/b}g)$
are alternate dual frames if and only if both
$h\perp E_{n/a}T_{m/b}g$, for all $(m,n)\not= (0,0)$, and $<h,g> = ab$.
\endproclaim

We proceed with the corresponding result to Theorem 2.2 for  
alternate dual frames.  $(1)\Leftrightarrow (2)$ in the theorem 
below was first proved by Janssen \cite{6}, section
1.3.2.  

\proclaim{Theorem 4.1}
For $g,h\in${\bf PF} and $(E_{mb}T_{na}g)$ a WH-frame for $L^{2}(R)$,
 the following are equivalent:
\vskip5pt
(1)  $(E_{mb}T_{na}h)_{m,n\in Z}$ is an alternate dual frame for
$(E_{mb}T_{na}g)_{m,n\in Z}$.
\vskip5pt
(2)  We have:

\ \ \ \ \ (a)  $\sum_{n\in Z}h(x-na)\overline{g(x-na-k/b)} = 0$  
a.e. for all
$k\not= 0$.
\vskip5pt
\ \ \ \ \ (b)  $\sum_{n\in Z}h(x-na)\overline{g(x-na)} = b$, a.e.
\vskip5pt
(3)  $h = S^{-1}g+f$, where $f\in L^{2}(R)$ and $f\perp  
\text{span}_{n,m\in Z}E_{n/a}T_{m/b}g$.
\endproclaim

\demo{Proof}
$(1)\Leftrightarrow (2)$:  This is Proposition 2.1 combined with the
theorem of Wexler-Raz.

$(1)\Rightarrow (3)$:  By the Wexler-Raz Theorem, $h\perp  
E_{n/a}T_{m/b}g$,
for all $(m,n)\not= (0,0)$.  By Corollary 2.3, we also have that
$S^{-1}g\perp E_{n/a}T_{m/b}g$.  Hence, $f=h-S^{-1}g\perp   
E_{n/a}T_{m/b}g$.  Again by the Wexler-Raz Theorem, $<h,g> = ab$ and  
applying Corollary 2.3 again,
$$
<h-S^{-1}g,g> = <h,g>+<S^{-1}g,g>  = ab-ab = 0.
$$
It follows that $h = S^{-1}g + (h-S^{-1}g) = S^{-1}g + f$ and  
$f\perp E_{n/a}T_{m/b}g$, for
all $n,m\in Z$.

$(3)\Rightarrow (1)$:  Fix $(m,n)\not= (0,0)$.  We compute using
Corollary 2.3:
$$
<h,E_{n/a}T_{m/b}g> = <S^{-1}g+f,E_{n/a}T_{m/b}g> =
$$
$$
<S^{-1}g,E_{n/a}T_{m/b}g> +
<f,E_{n/a}T_{m/b}g> = 0+0=0.
$$
Also, using Corollary 2.3,
$$
<h,g> = <S^{-1}g+f,g> = <S^{-1}g,g>+<f,g> = ab + 0 = ab.
$$
So this implication follows from the Wexler-Raz Theorem.
\enddemo

Note that $(E_{mb}T_{na}g)_{m,n\in Z}$ is a normalized tight frame  
if and
only if we can replace $S^{-1}g$ in Theorem 4.1 by the function $g$ and in this
case Theorem 4.1 reduces to Theorem 2.2.  Also in this case, $S=I$ so
part (3) of the theorem becomes:  $h=g+f$ where $f\perp
\text{span}_{n,m\in Z}E_{n/a}T_{m/b}g$.

\proclaim{ACKNOWLEDGEMENT}
The authors express their deepest gatitude to A.J.E.M. Janssen for
making extensive recommendations for improvements of this manuscript.
\endproclaim

\enddocument